\documentclass[a4paper,12pt]{amsart}
\usepackage{amsmath}
\usepackage{amsthm}
\usepackage{amscd}
\usepackage{amssymb}
\usepackage{amsfonts}
\usepackage{latexsym}
\usepackage[mathscr]{eucal}

\newtheorem{thm}{Theorem}[section]
\newtheorem{prop}[thm]{Proposition}
\newtheorem{lemma}[thm]{Lemma}
\newtheorem{cor}[thm]{Corollary}

\def\supp{\mathop{\rm {supp}}\nolimits}
\def\graph{\mathop{\rm {graph}}\nolimits}
\def\length{\mathop{\rm {length}}\nolimits}

\begin{document}
\title{$C^*$-algebras associated with complex dynamical systems}
\author{Tsuyoshi Kajiwara}
\address[Tsuyoshi Kajiwara]{Department of Environmental and 
Mathematical Sciences, 
Okayama University, Tsushima, 700-8530,  Japan}      

\author{Yasuo Watatani}
\address[Yasuo Watatani]{Department of Mathematical Sciences, 
Kyushu University, Hakozaki, 
Fukuoka, 812-8581,  Japan}
\maketitle
\begin{abstract}
Iteration of a rational function $R$ gives a complex dynamical 
system on the Riemann sphere.  We introduce a $C^*$-algebra 
${\mathcal O}_R$ associated with $R$ as a Cuntz-Pimsner 
algebra of a Hilbert bimodule over the algebra $A = C(J_R)$
of continuous functions on the Julia set $J_R$ of $R$. The 
algebra ${\mathcal O}_R$ is a certain analog of the crossed 
product by a boundary action. We show that if the degree of 
$R$ is at least two, then $C^*$-algebra ${\mathcal O}_R$ 
is simple and purely infinite. For example if $R(z) = z^2 - 2$, 
then the Julia set $J_R = [-2,2]$ and the restriction  
$R : J_R \rightarrow J_R$ is topologically conjugate  to the 
tent map on $[0,1]$. The algebra ${\mathcal O}_{z^2 - 2}$ is 
isomorphic to the Cuntz algebra ${\mathcal O}_{\infty}$.  
We also show that the Lyubich measure associated with $R$
gives a unique KMS state on the $C^*$-algebra ${\mathcal O}_R$ 
for the gauge action at inverse temperature $\log (\deg R)$ if the 
Julia set contains no critical points. 
\end{abstract}

\section{Introduction}
For a branched covering $\pi : M \rightarrow M$, 
Deaconu and Muhly \cite{DM} introduced a $C^*$-algebra 
$C^*(M,\pi)$ as the $C^*$-algebra of the r-discrete 
groupoid constructed by Renault \cite{R}. In particular 
they consider rational functions  on the Riemann sphere 
$\hat{\mathbb C}$ and compute the K-groups of the 
$C^*$-algebra. See also the previous work \cite{D}
by Deaconu on $C^*$-algebras associated with continuous 
graphs and Anantharaman-Delaroche's  work \cite{A} on purely 
infinite $C^*$-algebras for expansive dynamical systems. 
Although  Deaconu and Muhly's work itself is interesting, 
we introduce a slightly different 
$C^*$-algebras 
${\mathcal O}_R(\hat{\mathbb C})$,  
${\mathcal O}_R = {\mathcal O}_R(J_R)$ and ${\mathcal O}_R(F_R)$
associated with a rational function $R$ on the Riemann sphere,  
the Julia set $J_R$ and the Fatou set $F_R$ of $R$ in this note.
The $C^*$-algebra
${\mathcal O}_R(\hat{\mathbb C})$ is defined 
as a Cuntz-Pimsner algebra of the Hilbert bimodule 
$Y= C(\graph R)$ over $B = C(\hat{\mathbb C})$. 
The $C^*$-algebra ${\mathcal O}_R$ is defined  
as a Cuntz-Pimsner algebra of the Hilbert bimodule $
X= C(\graph R|_{J_R})$ over $A = C(J_R)$. And the $C^*$-algebra 
${\mathcal O}_R(F_R)$ is similarly defined. 
The difference between Deaconu and Muhly's construction and ours 
is the following: 
They exclude branched points to construct 
their groupoids.  
We include branched points to construct our bimodules. 
An advantage of theirs
is that their  $C^*$-algebra is constructed from both a groupoid and 
a bimodule.  An advantage of ours is that our $C^*$-algebra 
${\mathcal O}_R$ is always simple and purely infinite if 
the degree of $R$ is at least two.

For example, if $R(z) = z^2 -2$, 
then the Julia set $J_R = [-2,2]$
and the restriction $R|_{J_R}: J_R \rightarrow  J_R$ is 
topologically conjugate to the tent map 
$h: [0,1] \rightarrow [0,1]$.  Then their $C^*$-algebra 
$C^*([0,1],h)$ is not simple and 
${\rm K}_0(C^*([0,1],h)) = {\mathbb Z} \oplus {\mathbb Z}$ and 
${\rm K}_1(C^*([0,1],h)) = \{0\}$. 
Our $C^*$-algebra ${\mathcal O}_{z^2-2}$ is simple and 
purely infinite and 
${\rm K}_0({\mathcal O}_{z^2-2}) = {\mathbb Z}$ and 
${\rm K}_1({\mathcal O}_{z^2-2}) = \{0\}$.  In fact the 
algebra ${\mathcal O}_{z^2-2}$ is isomorphic to 
the Cuntz algebra ${\mathcal O}_{\infty}$.

Even if $R$ is a quadratic polynomial $R(z) = z^2 -c$, 
the structure of the $C^*$-algebra ${\mathcal O}_R$ is 
closely related to the  property of $R$ as the complex dynamical 
system. If $c$ is not in the Mandelbrot set ${\mathcal M}$, 
then $C^*$-algebra ${\mathcal O}_R$ is isomorphic to the 
Cuntz algebra ${\mathcal O}_2$. 
If $c$ is in the interior of the main cardioid, then  
$C^*$-algebra ${\mathcal O}_R$ is not isomorphic to 
${\mathcal O}_2$.  In fact we have 
$K_0({\mathcal O}_R) = {\mathbb Z}$ and 
$K_1({\mathcal O}_R) = {\mathbb Z}$. 

We compare our construction with other important constructions  
by C. Delaroche \cite{A} and M. Laca - J. Spielberg \cite{LS}.
 They showed that 
a certain boundary action of a Kleinian group 
on the limit set yields a simple nuclear purely infinite 
$C^*$-algebra as groupoid $C^*$-algebra or crossed product.  
Recall that Sullivan's dictionary says that there is 
a strong analogy between the limit set $\Lambda_{\Gamma}$ of a 
Kleinian group $\Gamma$ 
and the Julia set $J_R$ of a rational function $R$.
The algebra ${\mathcal O}_R$ is 
generated by $C(J_R)$ and $\{S_f ; f \in C(\graph R|_{J_R})\}$. 
We regard the algebra ${\mathcal O}_R$ as 
a certain analog of the 
crossed product $C(\Lambda_{\Gamma}) \rtimes \Gamma$ of 
$C(\Lambda_{\Gamma})$ by a boundary action of a Klein group $\Gamma$. 
In fact the crossed product is generated by  $C(\Lambda_{\Gamma})$
and $\{\lambda _g : g \in \Gamma \}$.  Moreover commutation 
relations $\alpha (a)S_f = S_fa$ for $a \in C(J_R)$ and 
$\alpha _g (a) \lambda _g = \lambda _g a$ for $a \in C(\Lambda_{\Gamma})$
 are similar, where $\alpha (a)(x) = a(R(x))$ and 
$\alpha _g (a)(x) = a(g^{-1}x)$.  

Applying  a result by Fowler, Muhly and Raeburn \cite{FMR}, 
the quotient algebra 
${\mathcal O}_R(\hat{\mathbb C})/{\mathcal I} (I)$ by the 
ideal ${\mathcal I} (I)$ corresponding to the Fatou set is 
canonically isomorphic to ${\mathcal O}_R = {\mathcal O}_R(J_R)$. 
 
Several people \cite{DPZ}, \cite{MS}, \cite{S}, \cite{KPW1} 
considered conditions for simplicity of Cuntz-Pimsner algebras.  
We directly show that  $C^*$-algebra ${\mathcal O}_R$ is simple 
purely infinite through analyzing  branched coverings. 
A criterion by Schweizer \cite{S} can be applied to show 
the only simplicity of $C^*$-algebra ${\mathcal O}_R$. 
$C^*$-algebra ${\mathcal O}_R$ is separable and nuclear, 
and belongs to the UCT class.  Thus if $\deg R \geq 2$, 
then the isomorphisms class of ${\mathcal O}_R$ is completely 
determined by the $K$-theory together with the class of the 
unit by the classification theorem by Kirchberg-Phillips \cite{Ki},
\cite{Ph}.
If there exist no critical points in $J_R$, then 
simplicity of ${\mathcal O}_R$ is also given by the simplicity 
of a certain crossed product by 
an endomorhism with a transfer operator in 
Exel-Vershik \cite{EV}, \cite{E}.
Katsura \cite{Ka} studies his new construction which contains 
${\mathcal O}_R$ if there exist no critical points in $J_R$.
A difficulty of our analysis stems from the fact that 
the Julia set $J_R$ contains the critical points (i.e. 
the branched points. )

Many examples of Cuntz-Pimsner algebras arise from Hilbert 
bimodules which are finitely generated and projective as 
right module and the image of the left action is contained 
in the compacts. 
Another extremal case is studied by A. Kumjian in \cite{Ku}, 
where the image of the left action has trivial intersection 
with the compacts. 
Our case is a third one and the intersection of the image of 
the left action and the compacts is exactly represented  by 
the branched points, i.e., critical points of $R$.  Hence our 
study is focused on the behaviour of the branched points 
under iteration.  

We also showed that the Lyubich measure associated with $R$ 
gives a unique KMS state on 
the $C^*$-algebra ${\mathcal O}_R$ for the gauge action at 
inverse temperature $\log \deg R$ if the Julia set contains 
no critical points.  But the problem of KMS state is subtle 
and difficult, if the Julia set contains critical points.  We 
will discuss the case elsewhere.

We should note that B. Brenken also announced a study 
of Cuntz-Pimsner algebras associated with branched 
coverings in a satellite conference at Chenge, China in 2002, 
when we announced a content of the paper in the same conference. 

\section{Rational functions and Hilbert bimodules}
   We recall some facts on iteration of rational functions. 
Let $R$ be a rational function of the form 
$R(z) = \frac{P(z)}{Q(z)}$ with relative prime polynomials 
$P$ and $Q$.   The degree of $R$ is denoted by 
$d = \deg R := \max \{ \deg P, \deg Q \}$. 
If $\deg R = 1$ and the Julia set is not empty, 
then it is, in fact,  one point.  Therefore 
the classical Toeplitz algebra appears and the  
$C^*$-algebra ${\mathcal O}_R$ becomes $C({\mathbb T})$. 
Hence  we need to assume that $\deg R \geq 2$ to consider 
non-trivial ones. 
We regard a rational function $R$ as a $d$-fold branched 
covering map  $R : \hat{\mathbb C} \rightarrow \hat{\mathbb C}$
on the Riemann sphere $\hat{\mathbb C} = {\mathbb C} 
\cup \{ \infty \}$.  The sequence $(R^n)_n$ of iterations of $R$ 
gives a complex analytic dynamical system on $\hat{\mathbb C}$. 
The Fatou set $F_R$ of $R$ is the maximal open subset of 
$\hat{\mathbb C}$ on which $(R^n)_n$ is equicontinuous (or 
a normal family), and the Julia set $J_R$ of $R$ is the 
complement of the Fatou set in $\hat{\mathbb C}$.  

We always assume that a rational function $R$ is not a constant 
function.  Recall that a {\it critical point} of $R$  
is a point   
$z_0$  at which $R$ is not locally one to one.  
It is a zero of $R'$ or 
a pole of $R$ of order two or higher.  The image $w_0 =R(z_0)$ 
is called a {\it critical value} of $R$.    
Using appropriate local charts, if $R(z) = w_0 + c(z - z_0)^n + 
(\text{higher terms})$ with 
$n \geq 1$ and $c \not= 0$ on some neighborhood of $z_0$, 
then the integer $n = e(z_0) = e_R(z_0)$ is called the 
{\it branch index} of $R$ at $z_0$.  Thus $e(z_0) \geq 2$ if 
$z_0$ is a critical point and $e(z_0) = 1$ otherwise. Therefore 
$R$ is an $e(z_0) :1$ map in a punctured neighborhood of $z_0$.
By the Riemann-Hurwitz formula, there exist $2d - 2$ critical points 
counted with multiplicity, that is, 
$$
\sum _{z \in \hat{\mathbb C}} e(z) -1 = 2 \deg R -2 .
$$
Furthermore for each $w \in \hat{\mathbb C}$, we have 
$$
\sum _{z \in R^{-1}(w)} e(z) = \deg R .
$$ 
Let $C$ be the 
set of critical points of $R$ and $R(C)$ be the 
set of the critical values of $R$. We put  
$B = R^{-1}(R(C))$. Then the restriction 
$R : \hat{\mathbb C} \setminus B \rightarrow 
\hat{\mathbb C} \setminus R(C) $ is a $d :1$ regular covering, 
where $d = \deg R$. This means that any point $y \in 
\hat{\mathbb C} \setminus R(C) $ has an open neighborhood $V$ 
such that $R^{-1}(V)$ has $d$ connected components $U_1, \dots, 
U_d$ and the restriction $R|_{U_k} : U_k \rightarrow V$ is a homeomorphism 
for $k = 1, \dots , d$.  Thus $R$ has $d$ analytic local cross 
sections $S_k = (R|_{U_k})^{-1}$.  

We recall Cuntz-Pimsner algebras \cite{Pi}.  
Let $A$ be a $C^*$-algebra 
and $X$ be a Hilbert right $A$-module.  We denote by $L(X)$ be 
the algebra of the adjointable bounded operators on $X$.  For 
$\xi$, $\eta \in X$, the "rank one" operator $\theta _{\xi,\eta}$
is defined by $\theta _{\xi,\eta}(\zeta) = \xi(\eta|\zeta)$
for $\zeta \in X$. The closure of the linear span of rank one 
operators is denoted by $K(X)$.   We call that 
$X$ is a Hilbert bimodule over $A$ if $X$ is a Hilbert right  $A$-
module with a homomorphism $\phi : A \rightarrow L(X)$.  We assume 
that $X$ is full and $\phi$ is injective.

   Let $F(X) = \oplus _{n=0}^{\infty} X^{\otimes n}$
be the Fock module of $X$ with a convention 
$X^{\otimes 0} = A$. 
 For $\xi \in X$, the creation operator
$T_{\xi} \in L(F(X))$ is defined by 
$$
T_{\xi}(a) =  \xi a  \qquad \text{and } \ 
T_{\xi}(\xi _1 \otimes \dots \otimes \xi _n) = \xi \otimes 
\xi _1 \otimes \dots \otimes \xi _n .
$$
We define $i_{F(X)}: A \rightarrow L(F(X))$ by 
$$
i_{F(X)}(a)(b) = ab \qquad \text{and } \ 
i_{F(X)}(a)(\xi _1 \otimes \dots \otimes \xi _n) = \phi (a)
\xi _1 \otimes \dots \otimes \xi _n 
$$
for $a,b \in A$.  The Cuntz-Toeplitz algebra ${\mathcal T}_X$ 
is the $C^*$-algebra on $F(X)$ generated by $i_{F(X)}(a)$
with $a \in A$ and $T_{\xi}$ with $\xi \in X$.  
Let $j_K : K(X) \rightarrow {\mathcal T}_X$ be the homomorphism 
defined by $j_K(\theta _{\xi,\eta}) = T_{\xi}T_{\eta}^*$. 
We consider the ideal $I_X := \phi ^{-1}(K(X))$ of $A$. 
Let ${\mathcal J}_X$ be the ideal of ${\mathcal T}_X$ generated 
by $\{ i_{F(X)}(a) - (j_K \circ \phi)(a) ; a \in I_X\}$.  Then 
the Cuntz-Pimsner algebra ${\mathcal O}_X$ is the 
the quotient ${\mathcal T}_X/{\mathcal J}_X$ . 
Let $\pi : {\mathcal T}_X \rightarrow {\mathcal O}_X$ be the 
quotient map.  Put $S_{\xi} = \pi (T_{\xi})$ and 
$i(a) = \pi (i_{F(X)}(a))$. Let
$i_K : K(X) \rightarrow {\mathcal O}_X$ be the homomorphism 
defined by $i_K(\theta _{\xi,\eta}) = S_{\xi}S_{\eta}^*$. Then 
$\pi((j_K \circ \phi)(a)) = (i_K \circ \phi)(a)$ for $a \in I_X$.   
We note that  the Cuntz-Pimsner algebra ${\mathcal O}_X$ is 
the universal $C^*$-algebra generated by $i(a)$ with $a \in A$ and 
$S_{\xi}$ with $\xi \in X$  satisfying that 
$i(a)S_{\xi} = S_{\phi (a)\xi}$, $S_{\xi}i(a) = S_{\xi a}$, 
$S_{\xi}^*S_{\eta} = i((\xi | \eta)_A)$ for $a \in A$, 
$\xi, \eta \in X$ and $i(a) = (i_K \circ \phi)(a)$ for $a \in I_X$.
We usually identify $i(a)$ with $a$ in $A$.  We denote by 
${\mathcal O}_X^{alg}$ the $\ ^*$-algebra generated algebraically 
by $A$  and $S_{\xi}$ with $\xi \in X$. There exists an action 
$\gamma : {\mathbb R} \rightarrow Aut \ {\mathcal O}_X$
with $\gamma_t(S_{\xi}) = e^{it}S_{\xi}$, which is called the  
gauge action. Since we assume that $\phi: A \rightarrow L(X)$ is 
isometric, there is an embedding $\phi _n : L(X^{\otimes n})
 \rightarrow L(X^{\otimes n+1})$ with $\phi _n(T) = 
T \otimes id_X$ for $T \in L(X^{\otimes n})$ with the convention 
$\phi _0 = \phi : A \rightarrow L(X)$.  We denote by ${\mathcal F}_X$
the $C^*$-algebra generated by all $K(X^{\otimes n})$, $n \geq 0$ 
in the inductive limit algebra $\varinjlim L(X^{\otimes n})$. 
Let ${\mathcal F}_n$ be the $C^*$-subalgebra of ${\mathcal F}_X$ generated by 
$K(X^{\otimes k})$, $k = 0,1,\dots, n$, with the convention 
${\mathcal F}_0 = A = K(X^{\otimes 0})$.  Then  ${\mathcal F}_X = 
\varinjlim {\mathcal F}_n$.

Let $\graph R = \{(x,y) \in \hat{\mathbb C}^2 ; y = R(x)\} $ 
be the graph of a rational function $R$.  Consider a $C^*$-algebra 
$B = C(\hat{\mathbb C})$. Let $Y = C(\graph R)$.  
Then $Y$ is a $B$-$B$ bimodule by 
$$
(a\cdot f \cdot b)(x,y) = a(x)f(x,y)b(y)
$$
for $a,b \in B$ and $f \in Y$. We introduce a $B$-valued 
inner product $(\ |\ )_B$ on $Y$ by 
$$
(f|g)_B(y) = \sum _{x \in R^{-1}(y)} e(x) \overline{f(x,y)}g(x,y)
$$
for $f,g \in Y$ and $y \in \hat{\mathbb C}$.  We need branch index 
$e(x)$ in the formula of the inner product above.
Put $\|f\|_2 = \|(f|f)_B\|_{\infty}$. 

\begin{lemma}
   The above $B$-valued inner product is well defined, 
that is, $\hat{\mathbb C} \ni y \mapsto (f|g)_B(y) 
\in {\mathbb C}$ is continuous.  
\label{lemma:1}
\end{lemma}  

\begin{proof} 
Let $d = \deg R$.  If $y_0$ is not a 
critical value of $R$, then there are $d$ distinct  
continuous cross sections, say $S_1, \dots , S_d$, of $R$
defined on some open neighbourhood $V$ of $y_0$ 
such that $V$ contains no critical values of $R$. 
Then for $y \in V$, 
$$
(f|g)_B(y) = \sum _{i = 1}^d  \overline{f(S_i(y),y)}g(S_i(y),y).
$$
Thus, $(f|g)_B(y)$ is continuous at $y_0$.  Next consider 
the case that  $y_0$ is a critical value of $R$.  Let 
$x_1, \dots , x_r$ be the distinct points of $R^{-1}(y_0)$ 
Then there exist open neighbourhoods $V$ of $y_0$ and $U_j$ 
of $x_j$ such that $R : U_j \rightarrow V$ is a $e(x_j):1$ 
branched covering and is expressed as $f(z) = z^{e(x_j)}$ 
by local charts for $j = 1, \dots , r$.  Consider 
$S_k^{(j)}(z) = r^{1/e(x_j)}\exp (i(\frac{\theta}{e(x_j)} + 
\frac{2k\pi}{e(x_j)}))$ for $z = r\exp (i\theta)$ 
and $S_k^{(j)}(0) = 0$ under the local coordinate. 
Thus there exist cross sections, say $S_k^{(j)}$  of $R$ 
$k = 1, \dots , e(x_j)$  for each $j = 1, \dots , r$ 
defined on $U$ of $y_0$  such any $S_k^{(j)}$ is 
continuous at least at $y_0$ but not necessarily continuous 
on $V$.  
Then for $y \in V$, 

\begin{align*}
\lim _{y \rightarrow y_0} (f|g)_B(y) 
& = \lim _{y \rightarrow y_0} 
\sum _{j = 1}^r  \sum _{k = 1}^{e(x_j)} 
 \overline{f(S_k^{(j)}(y),y)}g(S_k^{(j)}(y),y) \\
& = \sum _{j = 1}^r  \sum _{k = 1}^{e(x_j)} 
 \overline{f(S_k^{(j)}(y_0),y_0)}g(S_k^{(j)}(y_0),y_0) \\
& = \sum _{j = 1}^r  {e(x_j)} 
 \overline{f(x_j,y_0)}g(x_j,y_0) 
= (f|g)_B(y_0)
\end{align*}

Thus $(f|g)_B(y)$ is continuous at $y_0$.

\end{proof} 
  
The left multiplication of $B$ on $Y$ gives 
the left action $\phi : B \rightarrow L(Y)$ 
such that 
$(\phi (b)f)(x,y) = b(x)f(x,y)$ for $b \in B$ and $f \in Y$.

\begin{prop}If $R$ is a rational function, then 
$Y = C(\graph R)$ is a full Hilbert bimodule over 
$B = C(\hat{\mathbb C})$ without completion. 
The left action $\phi : B \rightarrow L(Y)$ is unital and 
faithful. 
\end{prop}

\begin{proof} Let $d = \deg R$.  
For any $f \in Y = C(\graph R)$, we have 
$$
\| f\|_{\infty} \leq \| f \|_2 
= (\sup _y  \sum _{x \in R^{-1}(y)} e(x)|f(x,y)|^2)^{1/2}
\leq \sqrt{d} \| f\|_{\infty}
$$
Therefore two norms $\|\ \|_2$ and $\|\ \|_{\infty}$ are 
equivalent.  Since $C(\graph R)$ is complete  with  respect to 
$\| \ \|_{\infty}$, it is also complete with  respect to 
$\|\ \|_2$.

Since $(1|1)_B(y) =  \sum _{x \in R^{-1}(y)} e(x)1 = \deg R$, 
$(Y|Y)_B$ contains the identity $I_B$ of $B$.  Therefore 
$Y$ is full. If $b \in B$ is not zero, then there exists 
$x_0 \in \hat{\mathbb C}$ with $b(x_0) \not = 0$.  Choose 
$f \in Y$ with $f(x_0,R(x_0)) \not= 0$.  Then 
$\phi (b)f \not= 0$.  Thus $\phi$ is faithful.     

\end{proof} 

Since the Julia set $J_R$ is completely invariant, i.e., 
$R(J_R) = J_R = R^{-1}(J_R)$, we can consider the restriction 
$R|_{J_R} : J_R \rightarrow J_R$, which will be often denoted by 
the same letter $R$.
 Let $\graph R|_{J_R} = \{(x,y) \in J_R \times J_R \ ; \ y = R(x)\} $
be the graph of the restricton map $R|_{J_R}$. 
Let $A = C(J_R)$ and $X = C(\graph R|_{J_R})$. Through 
a restriction of the above action,  $X$ is an $A$-$A$ bimodule. 

\begin{cor}
Let $R$ be a rational function with $J_R \not= \phi$, for example, 
$\deg R \geq 2$.  
Then $X = C(\graph R|_{J_R})$ is a full Hilbert bimodule over 
$A = C(J_R)$ without completion, where $A$-valued inner product 
$(\ |\ )_A$ on $X$ is given by 
$$
(f|g)_A(y) = \sum _{x \in R^{-1}(y)} e(x) \overline{f(x,y)}g(x,y)
$$
for $f,g \in X$ and $y \in J_R$.  
The left action $\phi : A \rightarrow L(X)$ is unital and 
faithful.
\end{cor} 

\begin{proof}
It is an immediate consequence of the fact that 
$R(J_R) = J_R = R^{-1}(J_R)$. 
\end{proof} 

\noindent
{\bf Definition.} 
We introduce the $C^*$-algebras 
${\mathcal O}_R(\hat{\mathbb C})$, 
${\mathcal O}_R(J_R)$ and ${\mathcal O}_R(F_R)$
associated with a rational function $R$. 
The $C^*$-algebra
${\mathcal O}_R(\hat{\mathbb C})$ is defined 
as a Cuntz-Pimsner algebra of the Hilbert bimodule 
$Y= C(\graph R)$ over 
$B = C(\hat{\mathbb C})$. 
When the Julia set $J_R$ is not empty, for example 
$\deg R \geq 2$, we introduce
the $C^*$-algebra ${\mathcal O}_R(J_R)$ 
as a Cuntz-Pimsner algebra of the Hilbert bimodule $
X= C(\graph R|_{J_R})$ over $A = C(J_R)$. 
When the Fatou set $F_R$ is not empty, the $C^*$-algebra
${\mathcal O}_R(F_R)$ is defined similarly.
Sullivan's dictionary says that there is 
a strong analogy between the limit set of a Kleinian group 
and the Julia set of a rational function.  Therefore 
we simply denote by ${\mathcal O}_R$ the $C^*$-algebra 
${\mathcal O}_R(J_R)$ to emphasize the analogy.  
 
\begin{prop}
Let $R$ be a rational function and  
$Y = C(\graph R)$ be a Hilbert bimodule over 
$B = C(\hat{\mathbb C})$.  Then there exists an 
isomorphism 
$$
\varphi : Y^{\otimes n} \rightarrow C(\graph R^n) 
$$
as a Hilbert bimodule over $B$ such that 
\begin{align*}
& (\varphi (f_1 \otimes \dots \otimes f_n))(x,R^n(x)) \\ 
& = f_1(x,R(x))f_2(R(x),R^2(x))\dots f_n(R^{n-1}(x),R^n(x))
\end{align*}
for $f_1,\dots, f_n \in Y$ and $x \in \hat{\mathbb C}$.  
Moreover when the Julia set $J_R$ is not empty, 
let $A = C(J_R)$ and $X = C(\graph R|_{J_R})$.
Then we have a similar isomorphism 
$$
\varphi : X^{\otimes n} \rightarrow C(\graph R^n|_{J_R}) 
$$
as a Hilbert bimodule over $A$. Here we use the same symbol
$\varphi$ to save the notation. 
\end{prop}

\begin{proof}
It is easy to see that $\varphi$ is well-defined
and a bimodule homomorphism. We show that $\varphi$ preserves
inner product. The point is to use the following chain rule 
for branch index: 
$$
e_{R^n}(x) = e_R(x)e_R(R(x))\dots e_R(R^{n-1}(x))
$$
Consider the case when $n = 2$ for simplicity of the notation.
\begin{align*}
& (f_1 \otimes f_2 | g_1 \otimes g_2)_B(x_2) 
 = (f_2|(f_1|g_1)_Bg_2)_B(x_2) \\
& = \sum _{x_1 \in R^{-1}(x_2)} e_R(x_1)
 \overline{f_2(x_1,x_2)}(f_1|g_1)_B(x_1)g_2(x_1,x_2) \\
& = \sum _{x_1 \in R^{-1}(x_2)} e_R(x_1)
 \overline{f_2(x_1,x_2)}
(\sum _{x_0 \in R^{-1}(x_1)} e_R(x_0)
\overline{f_1(x_0,x_1)}g_1(x_0,x_1)
g_2(x_1,x_2)) \\
& = \sum _{x_0 \in R^{-2}(x_2)} e_{R^2}(x_0)
\overline{f_1(x_0,R(x_0))f_2(R(x_0),x_2)}
g_1(x_0,R(x_0))g_2(R(x_0),x_2) \\
& = \sum _{x_0 \in R^{-2}(x_2)} e_{R^2}(x_0)
\overline{(\varphi (f_1 \otimes f_2))(x_0,x_2)} 
(\varphi (g_1 \otimes g_2))(x_0,x_2) \\
& = (\varphi (f_1 \otimes f_2) | \varphi (g_1 \otimes g_2))(x_2) 
\end{align*}  

Since $\varphi$ preserves inner product, $\varphi$ is one to one. 
The non-trivial one  is to show that $\varphi$ is onto. Since 
the image of $\varphi$ is a $\ ^*$-subalgebra of $C(\hat{\mathbb C})$
and separates the two points, the image of $\varphi$ is dense in 
$C(\hat{\mathbb C})$ with respect to $\| \ \|_{\infty}$ by the 
Stone-Wierstrass Theorem.  
Since two norms $\|\ \|_2$ and $\|\ \|_{\infty}$ are equivalent 
and $\varphi$ is isometric with respect to $\|\ \|_2$, 
$\varphi$ is onto.

We recall that the Julia set $J|_R$ of $R$ is completely invariant 
under $R$ and the Julia set of $R$ coincides with the Julia set of 
$R^n$. Therefore the proof is valid for $A$ and $X$ also.
\end{proof}  

The set $C$ of critical points of $R$  is described by the 
ideal  $I_Y = \phi ^{-1}(K(Y))$ of $B$.

\begin{prop}
Let $R$ be a rational function and  
$Y = C(\graph R)$ be a Hilbert bimodule over 
$B = C(\hat{\mathbb C})$. 
Similarly let $A = C(J_R)$ and $X = C(\graph R|_{J_R})$.
Consider the ideal $I_Y = \phi ^{-1}(K(Y))$ of $B$  and 
the ideal $I_X = \phi ^{-1}(K(X))$ of $A$. Then  
$I_Y = \{b \in B ; b \ \text{vanishes on } C \}$.
and  $I_X = \{a \in A ; a \ \text{vanishes on } C \cap J_R \}$. 

\label{prop:critical}
\end{prop}

\begin{proof}Let $d = \deg R$, then there exist $2d-2$ critical 
points counted with multiplicity.  Therefore $C$ is a finite 
non-empty set. 
Firstly, let us take  $b \in B$ with a compact support $S = \supp (b)$ in 
$\hat{\mathbb C} \setminus C$. For any $x \in S$, since 
$x$ is not a critical point, there exists
an open neighbourhood $U_x$ of $x$ such that $U_x \cap C = \phi$ 
and the restriction $R|_{U_x}: U_x \rightarrow R(U_x)$ 
is a homeomorphism. Since $S$ is compact, there exists a finite 
subset $\{x_1,\dots , x_m\}$ such that 
$S \subset \cup _{i=1}^m U_{x_i} \subset \hat{\mathbb C} \setminus C$. 
Let $(f_i)_i$ be a finite family in $C(\hat{\mathbb C})$ such that 
$0 \leq f_i \leq 1$, $\supp(f_i) \subset U_{x_i}$ for 
$i=1,\dots ,m$ and $\sum _{i=1}^{m} f_i(x) = 1$ 
for $x \in S$.  Define 
$\xi _i, \eta _i \in C(\graph R)$ by 
$\xi _i(x,R(x)) = b(x)\sqrt{f_i(x)}$ and 
$\eta _i(x,R(x)) = \sqrt{f_i(x)}$.  Consider 
$T:= \sum _{i = 1}^k \theta _{\xi _i, \eta _i} \in K(Y)$.
We shall show that $T = \phi (b)$.      
For any $\zeta \in C(\graph R)$, we have 
$(\phi (b)\zeta )(x,y) = b(x)\zeta (x,y)$ and 
\begin{align*}
& (T\zeta )(x,y) = \sum _i \xi _i(x,y) 
\sum _{x' \in R^{-1}(y)} e(x')
\overline{\eta _i (x',y)}\zeta (x',y) \\
& = \sum _i b(x)\sqrt{f_i(x)} 
\sum _{x' \in R^{-1}(y)} e(x')\sqrt{f_i(x')}\zeta (x',y).
\end{align*}
In the case when $b(x) = 0$, we have 
$$
(T\zeta )(x,y) = 0 = (\phi (b)\zeta )(x,y).
$$
In the case when $b(x) \not= 0$, we have $x \in \supp (b) = S 
\subset \cup _{i=1}^m U_{x_i}$. Hence $e(x) = 1$ and 
$x \in U_{x_i}$ for some $i$. Put $y = R(x)$.  
Then for any $x' \in R^{-1}(y)$ 
with $x' \not= x$, $f_i(x') = 0$, because $x' \in U_{x_i}^c$. 
Therefore we have 
$$
(T\zeta )(x,y) 
= \sum _{\{i;x \in U_{x_i}\}} b(x)f_i(x)\zeta (x,y) 
= b(x)\zeta (x,y) = (\phi (b)\zeta )(x,y).
$$
Thus $\phi (b) = T \in K(Y)$.  Now for a general $b
\in B$ which vanishes on the 
critical points $C$ of $R$, there exists a sequence 
$(b_n)_n$ in $B$ with 
compact supports $\supp (b_n) \subset  
\hat{\mathbb C} \setminus C$ such that 
$\|b - b_n\|_{\infty} \rightarrow 0$.  
Hence $\phi (b) \in K(Y)$, i.e., $b \in I_Y$.

Conversely let $b \in B$ and $b(c) \not= 0$ for some 
critical point $c$.  We may assume that $b(c) = 1$. 
Put $m = e(c) \geq 2$. 
We need to show that $\phi(b) \not\in K(Y)$.  
On the contrary suppose that $\phi(b) \in K(Y)$.  
Then for $\varepsilon = \frac{1}{5\sqrt{d}}$, there exists a finite subset 
$\{\xi _i, \eta _i \in X ; i = 1,\dots ,N\}$ 
such that $\| \phi(b) - 
\sum _{i =1}^N \theta _{\xi _i, \eta _i} \| < \varepsilon$.  
Choose a sequence $(x_n)_n$ in $\hat{\mathbb C} \setminus C$ 
such that $x_n \rightarrow c$.  Since $x_n$ is not a critical 
point, there exists an open neighbourhood $U_n$ of $x_n$ 
such that $U_x \cap C = \phi$ 
and the restriction $R|_{U_n}: U_n \rightarrow R(U_n)$ 
is a homeomorphism. There exists $\zeta _n \in Y$ with 
$\supp \zeta _n  \subset \{(x,R(x)) ; x \in U_n\}$, 
$\zeta _n(x_n,R(x_n)) = 1$ and 
$0 \leq \zeta _n \leq 1$.  Then 
$\|\zeta _n \|_2 \leq \sqrt{d}$. For any $x' \in R^{-1}(R(x_n))$ 
with $x' \not= x_n$, $\zeta_n(x',R(x_n)) = 0$,
Hence 
\begin{align*}
& |b(x_n) - \sum _{i=1}^{N} \xi _i(x_n,R(x_n))
\overline{\eta _i(x_n,R(x_n))}| \\
& = |b(x_n) - 
\sum _{i=1}^{N} \xi _i(x_n,R(x_n)) 
\sum _{x' \in R^{-1}(R(x_n))} e(x')
\overline{\eta _i (x',R(x_n))}(\zeta _n) (x',R(x_n))\\
& = |((\phi(b) - \sum _{i =1}^N \theta _{\xi _i, \eta _i})
  \zeta _n)(x_n,R(x_n))| \\
& \leq \|(\phi(b) - \sum _{i =1}^N \theta _{\xi _i, \eta _i})
  \zeta _n\|_2 
\leq \|(\phi(b) - \sum _{i =1}^N \theta _{\xi _i, \eta _i}\|
  \|\zeta _n\|_2 
\leq \varepsilon \sqrt{d}.
\end{align*}  

Taking $n \rightarrow \infty$, we have 
$$
|b(c) - \sum _{i=1}^{N} \xi _i(c,R(c))
\overline{\eta _i(c,R(c))}| 
\leq \varepsilon \sqrt{d}.
$$
On the other hand, consider  $\zeta \in Y$ satisfying 
$\zeta (c,R(c)) = 1$,  $0 \leq \zeta \leq 1$ and 
$\zeta (x',R(x')) = 0$ 
for $x'\in R^{-1}(R(c))$ with $x' \not= c$
Then 
\begin{align*}
& |b(c) - \sum _{i=1}^{N} \xi _i(c,R(c))
e(c)\overline{\eta _i(c,R(c))}| \\
& = |b(c) - \sum _{i=1}^{N} \xi _i(c,R(c))
\sum _ {x' \in R^{-1}(R(c))}
e(x')\overline{\eta _i(x',R(x'))}\zeta (x',R(x'))| \\
& \leq \|(\phi(b) - \sum _{i =1}^N \theta _{\xi _i, \eta _i})
  \zeta \|_2 
\leq \varepsilon \sqrt{d}.
\end{align*}  

Since $e(c) \geq 2$ and $b(c) = 1$,  we have 
$$
\frac{1}{2} 
\leq |b(c) - \frac{1}{e(c)}b(c)| 
\leq \varepsilon \sqrt{d} + \frac{1}{e(c)}\varepsilon \sqrt{d} 
\leq 2\varepsilon \sqrt{d} = \frac{2}{5}
$$
This is a contradiction. 
Therefore  $\phi(b) \not\in K(Y)$
\end{proof}

\begin{cor}
$\ ^\#C = \dim (B/I_Y)$ and $\ ^\#(C \cap J_R) = \dim (A/I_X)$.
\end{cor}

\begin{cor}
The Julia set contains no critical points if and only if 
$\phi (A)$ is contained in $K(X)$ if and only if
$X$ is finitely generated projective right $A$ module.
\end{cor}

\section{Simplicity and pure infinteness}
Let $R$ be a rational function with the Julia set $J_R \not= \phi$. 
Let $A = C(J_R)$ and $X = C(\graph R|_{J_R})$. Define an 
endomorphism $\alpha : A \rightarrow A$ by 
$$
(\alpha (a))(x) = a(R(x))
$$
for $a \in A$, $x \in J_R$.  We also define a unital completely 
positive map $E_R: A \rightarrow A$ by
$$
(E_R(a))(y) := \frac{1}{\deg R} 
\sum _{x \in R^{-1}(y)} e(x)a(x)
$$
for $a \in A$, $y \in J_R$. In fact, for a constant function 
$\xi _0 \in X$ with 
$$
\xi _0 (x,y) := \frac{1}{\sqrt{\deg R}}
$$
we have     
$$
E_R(a) = (\xi _0|\phi (a)\xi _0)_A  
\text{and} \ E_R(I) = (\xi _0|\xi _0)_A = I. 
$$
We introduce an operator $D := S_{\xi _0} \in {\mathcal O}_R$.

\begin{lemma} 
In the above situation, for  $a, b  \in A$, we have the following:
\begin{enumerate}
\item $E_R (\alpha (a)b) = aE(b)$ and in particular 
$E_R (\alpha (a)) = a$
\item $D^*D = I$.
\item $\alpha (a) D = Da$. 
\item $D^*aD = E_R(a)$.
\item $\phi (A)\xi _0 = X$.
\end{enumerate}
\end{lemma}
 
\begin{proof}
(1): By the definition of $E_R$, we have 
\[
E_R (\alpha (a)b)(y) = \frac{1}{\deg R} 
\sum _{x \in R^{-1}(y)} e(x)a(R(x))b(x) = a(y)E(b)(y). 
\]
(2): $D^*D = S_{\xi _0}^*S_{\xi _0} = (\xi _0|\xi _0)_A = I.$ \\
\noindent
(3): We have  $\alpha (a) D = S_{\phi (\alpha (a))\xi _0}$ 
and $Da = S_{\xi _0 a}$.  Since 
$$
(\phi (\alpha (a))\xi _0)(x,R(x)) = (\alpha (a))(x)\xi _0(x,R(x)) 
= a(R(x)) = (\xi _0 a)(x,R(x)), 
$$¡¡ 
we have $\alpha (a) D = Da$. \\
\noindent
(4): $D^*aD = S_{\xi _0}^*aS_{\xi _0} = (\xi _0|\phi (a)\xi _0)_A 
= E_R(a)$. \\ 
\noindent
(5)For any $f \in X$, define $\tilde{f} \in A$ by 
$\tilde{f}(x) = \sqrt{\deg R}f(x,R(x))$.  
Then $\phi (\tilde{f}) \xi _0 = f$.  Thus $\phi (A)\xi _0 = X$.

\end{proof}

R. Exel introduced an interesting construction of  a 
crossed product by 
an endomorhism with a transfer operator in \cite{E}.

\begin{prop} 
$C^*$-algebra ${\mathcal O}_R$ is isomorphic to a Exel's crossed product by 
an endomorhism with a transfer operator $E_R$.
\end{prop} 
 
\begin{proof}By the universality, the Toeplitz algebras are isomorphic. 
Since 
$\phi (A)S_{\xi _0} = \{S_x ; x \in X\}$ can be identified 
with $X$ by (5), $(a,k)$ is a redundancy if and only if $k = \phi (a)$.    
Therefore $C^*$-algebra ${\mathcal O}_R$ is isomorphic to his 
algebra. 
\end{proof}

\noindent
{\bf Remark}.  If there exist no critical  
points in $J_R$, then the simplicity of ${\mathcal O}_R$ is 
a consequence of the  simplicity of his algebra, which is  
proved in Exel-Vershik \cite{EV}.  But we need a further argument 
because of the existence of critical points in $J_R$.   

\begin{lemma} 
In the same situation, for  $a \in A$ and 
$f_1,\dots , f_n \in X$, 
we have the following: 
$$
\alpha ^n(a)S_{f_1}\dots S_{f_n} = S_{f_1}\dots S_{f_n}a .
$$
\label{lemma:a-f}
\end{lemma}

\begin{proof}
It is enough to show that 
$\alpha (a)S_f = S_fa$ for $f \in X$. 
We have  $\alpha (a) S_f = S_{\phi (\alpha (a))f}$ 
and $S_fa = S_{fa}$.  Since 

\begin{align*}
& (\phi (\alpha (a))f)(x,R(x)) = (\alpha (a))(x)f(x,R(x)) \\
& = a(R(x))f(x,R(x))  = (fa)(x,R(x)), 
\end{align*}
 
\noindent
we have $\alpha (a)S_f = S_fa$. \\
\end{proof}

\par
\begin{lemma}
Let $R$ be a rational function with $\deg R \geq 2$.  
For any non-zero positive element $a \in A$ and for any  
$\varepsilon > 0$ there exist $n \in \mathbb{N}$ and 
$f \in X^{\otimes n}$  with $(f|f)_A = I$ such that  
$$
 \|a\| -\varepsilon \le S_f^* a S_f \le \|a\|
$$
\label{lemma:epsilon}
\end{lemma}
\begin{proof} Let $x_0$ be a point in $J_R$ with 
$|a(x_0)| = \| a \|$.  For any  
$\varepsilon > 0$ there exist an open neighbourhood $U$ of $x_0$
in $J_R$ such that for any $x \in U$  we have 
$\|a\| -\varepsilon  \le a(x) \le \|a\| $. 
Choose  anothter open neighbourhood $V$ of $x_0$ in $J_R$ and a compact 
subset $K \subset J_R$ satisfying $V \subset K \subset U$.
Since $\deg R \geq 2$, there exists $n \in \mathbb{N}$ such that 
$R^n(V) =J_R$ by Beardon \cite{B} Theorem 4.2.5.
We identify $X^{\otimes n}$ with $C(\graph R^n|_{J_R})$ 
as in Lemma \ref{lemma:a-f}. 
Define closed subsets $F_1$ and $F_2$ of $J_R \times J_R$ by  

\begin{align*}
 F_1 & = \{(x,y)\in J_R \times J_R | y = R^n(x), x \in K\} \\
 F_2 & = \{(x,y) \in J_R \times J_R | y = R^n(x), x \in U^{c} \}
\end{align*}

Since $F_1 \cap F_2 = \phi $, there exists  
$g \in C(\graph R^n|_{J_R})$ such that 
$0 \le g(x) \le 1$  and 
$$
 g(x,y) =\left\{\begin{array}{cc}
  1,  & (x,y)\in F_1  \\
  0,  & (x,y)\in F_2
       \end{array}\right.
$$
\par
Since $R^n(V)=J_R$, for any $y \in J_R$ there exists 
$x_1 \in V$ such that $R^n(x_1)=y$, so that 
$(x_1,y) \in F_1$. Therefore  
\begin{align*}
 (g|g)_A(y) & = \sum_{x \in (R^n)^{-1}(y)}e_{R^n}(x) |g(x,y)|^2 \\
            & \ge e_{R^n}(x_1) |g(x_1,y)|^2  \ge 1 
           \end{align*}

Let $b:=(g|g)_A$.  Then $b(y) = (g|g)_A(y) \ge 1$.  Thus 
$b \in A$ is positive and invertible.
We put $f := g b^{-1/2} \in X^{\otimes n}$.  Then  
$$
 (f|f)_A   = (g b^{-1/2}| g b^{-1/2})_A 
          = b^{-1/2} (g|g)_A b^{-1/2} 
           = I.
$$

\par
For any $y \in J_R$ and $x \in (R^n)^{-1}(y)$, if $x \in U$, 
then $\|a\|-\varepsilon \le a(x)$, and  
if $x \in U^c$, then  $f(x,y)=g(x,y) b^{-1/2}(y) = 0$. 
Therefore   
\begin{align*}
 \|a \|-\varepsilon & = (\|a\|-\varepsilon) (f|f)_A(y) \\
                    & = (\|a\|-\varepsilon) \sum_{x \in (R^n)^{-1}(y)} 
                     e_{R^n}(x)  |f(x,y)|^2 \\
                    & \leq \sum_{x \in (R^n)^{-1}(y)}
                      e_{R^n}(x)a(x)|f(x,y)|^2 \\
                    & = (f|af)_A(y) = S_f^* a S_f (y). 
\end{align*}

We also have that 
\[
S_f^* a S_f = (f|af)_A \leq \| a \| (f|f)_A = \| a \| .
\] 

\end{proof}

\begin{lemma}
Let $R$ be a rational function with $\deg R \geq 2$.  
For any non-zero positive element $a \in A$ and for any  
$\varepsilon > 0$ with $0 < \varepsilon < \|a\|$, 
there exist $n \in \mathbb{N}$ and $u \in X^{\otimes n}$ 
such that  
$$
 \| u \| _2 \le (\|a\| - \varepsilon)^{-1/2} \qquad \text{and} 
\quad  S_u^* a S_u = I
\label{lemma:u-epsilon}
$$
\end{lemma}
\begin{proof} 
For any $a \in A$ and $\varepsilon > 0$ as above, we choose 
$f \in (X)^{\otimes n}$ as in Lemma \ref{lemma:epsilon}. 
Put $c=S_f^*aS_f$.  
Since $0 < \|a\|-\varepsilon \le c \le \|a\|$, $c$ is positive 
and invertible.  Let $u:=fc^{-1/2}$. Then 
$$
 S_u^* a S_u = (u|au)_A = (fc^{-1/2}| a fc^{-1/2})_A 
             = c^{-1/2}(f|af)_A c^{-1/2} = I.
$$
Since $\|a\|-\varepsilon \le c$, 
we have $c^{-1/2}\le (\|a\|-\varepsilon)^{-1/2}$.¡¡
Hence 
$$
 \|u\| _2 = \|f c^{-1/2}\| _2 \le \|c^{-1/2}\| _2 
\leq (\|a\| - \varepsilon)^{-1/2}.
$$
\end{proof}

\begin{lemma}
Let $R$ be a rational function with $\deg R \geq 2$. 
For any $m \in \mathbb{N}$, any $n \in \mathbb{N}$,  
any $T \in L(X^{\otimes n})$ and 
any  $\varepsilon > 0$,  there exists a positive element $a \in A$ 
such that 
$$
 \| \phi(\alpha^n(a))T\|^2 \ge \|T\|^2 -\varepsilon
$$
and $a \alpha^k (a) = 0$ for $k=1$, $\cdots$, $m$.
\label{lemma:T-epsilon}
\end{lemma}

\begin{proof}
For any $m \in \mathbb{N}$, let  
${\mathcal P}_m = \{x \in J_R\ | R^k(x) =x, 
\quad \text{for some } k = 1,\dots , m\}$. 
Since $R$ is a rational function, ${\mathcal P}_m$ 
is a finite set.

For any $n \in \mathbb{N}$,  any $T \in L(X^{\otimes n})$ and 
any  $\varepsilon > 0$, there exists $f \in X^{\otimes n}$
such that $ \|f\|_2=1$ and 
$\|T\|^2 \ge \|Tf\|_2^2 > \|T\|^2 -\varepsilon $. 
We still identify $X^{\otimes n}$ with $C(\graph R^n|_{J_R})$.  
Then there exists $y_0 \in J_R$ such that  
$$
 \|Tf\|_2^2 =  \sum_{x \in (R^n)^{-1}(y_0)}
  e_{R^n}(x) |(Tf)(x,y_0)|^2 > \|T\|^2 -\varepsilon.
$$
Since $y \mapsto (Tf|Tf)_A(y)$ is continuous and 
$$
 \|Tf\|_2^2 = \sup_{y \in J_R} \sum_{x \in (R^n)^{-1}(y)}
e_{R^n}(x)|(Tf)(x,y)|^2,
$$ 
there exists an open neighbourhood $U_{y_0}$ of $y_0$ such that 
for any $y \in U_{y_0}$
$$
   \sum_{x \in (R^n)^{-1}(y)}
  e_{R^n}(x) |(Tf)(x,y)|^2 > \|T\|^2 -\varepsilon .
$$
Since $R$ is a rational function with $\deg R \geq 2$, 
the Julia set $J_R$ is a perfect set and uncountable by 
Beardon \cite{B} Theorem 4.2.4. Therefore  
 $U_{y_0} \backslash{\mathcal P}_m \ne \phi$.  
Choose  $y_1 \in U_{y_0} \backslash{\mathcal P}_m$. 
Since  $R^k(y_1)\not=  y_1$ ($k=1,\dots,m$),  
 $y_1 \not\in (R^k)^{-1}(y_1)$. 
For $k=1,\dots,m$, choose an open neighbourhood 
$W_k$ of $y_1$ such that 
$W_k \subset U_{y_0}$ and 
 $(R^k)^{-1}(W_k) \cap W_k = \phi$. 
Put  $W=\cap_{k=1}^{m}W_k$. 
Choose $a \in A = C(J_R)$ such that 
$a(y_1) = 1$, $0 \le a \le 1$ and  $\supp \ a \subset W$.  
Since $(R^k)^{-1}(W) \cap W = \phi$, 
$$
 \supp(\alpha^k(a)) \cap \supp (a) = \phi \qquad (k=1,\dots,m). 
$$
It implies that $\alpha^k(a) a =0$ ($k=1,\dots,m$). 

Since $y_1 \in  U_{y_0}$,     
\begin{align*}
\| \phi(\alpha^n(a))Tf\| _2^2  
 & =  \sup_{y \in J_R} \sum_{x \in (R^n)^{-1}(y)}
       e_{R^n}(x)|(a(R^n(x))(Tf)(x,y)|^2  \\
 & \ge \sum_{x \in (R^n)^{-1}(y_1)}e_{R^n}(x)|a(y_1)
        (Tf)(x,y_1)|^2 \\
 & = \sum_{x \in (R^n)^{-1}(y_1)}e_{R^n}(x)|(Tf)(x,y_1)|^2 \\
 & > \|T\|^2 -\varepsilon. 
\end{align*} 

Therefore we have $\| \phi(\alpha^n(a))T\|^2 \ge \|T\|^2 -\varepsilon$.

\end{proof}

Let ${\mathcal F}_n$ be the $C^*$-subalgebra of ${\mathcal F}_X$ 
generated by $K(X^{\otimes k})$, $k = 0,1,\dots, n$ 
and $B_n$ be the $C^*$-subalgebra of ${\mathcal O}_X$ 
generated by 
\[
\bigcup_{k=1}^n \{S_{x_1} \dots S_{x_k}S_{y_k}^* \dots S_{y_1}^* : 
x_1, \dots x_k,  y_1, \dots y_k \in X \} \cup A.
\]
 In the following Lemma \ref{lemma:I-free} we shall use 
an isomorphism 
$\varphi : {\mathcal F}_n \rightarrow B_n$ 
as in Pimsner \cite{Pi}  and  Fowler-Muhly-Raeburn \cite{FMR}
such that 
$$
\varphi (\theta_{x_1 \otimes \dots \otimes x_k, 
         y_1 \otimes \dots \otimes y_k}) 
    =  S_{x_1} \dots S_{x_k}S_{y_k}^* \dots S_{y_1}^*.
$$
To simplify notation, we put $S_x = S_{x_1} \dots S_{x_k}$ 
for $x = x_1 \otimes \dots \otimes x_k \in X^{\otimes k}$ .

\begin{lemma}
Let $R$ be a rational function with $\deg R \geq 2$. 
Let $b = c^*c$ for some $c \in {\mathcal O}_X^{alg}$.  
We decompose  $b = \sum _j b_j$ with 
$\gamma _t(b_j) = e^{ijt}b_j$.
For any  $\varepsilon >0 $
there exists $P \in A$ with $0\le P \le I$ satisfying the 
following:  
\begin{enumerate}
 \item $Pb_jP = 0$ \qquad $(j\ne 0)$
 \item $\|Pb_0P\| \ge \|b_0\| -\varepsilon $
\end{enumerate} 
\label{lemma:I-free}
\end{lemma} 

\begin{proof} For $x \in X^{\otimes n}$, we define $\length (x) = n$ 
with the convention $\length (a) = 0$ for $a \in A$.  
We write $c$ as a finite sum $c = a + \sum _i S_{x_i}S_{y_i}^*$.
Put $m = 2 \max \{\length (x_i), \length (y_i) ; i\}$.  

\par\noindent
For $j > 0$, each $b_j$ is a finite sum of terms in the form such that 
$$
S_x S_y^* \qquad x \in X^{\otimes (k+j)}, \qquad y \in X^{\otimes k} 
\qquad 0 \le k+j \le m 
$$
In the case when $j<0$, 
$b_j$ is a finite sum of terms in the form such that 
$$
 S_x S_y^* \qquad x \in X^{\otimes k}, \qquad y \in X^{\otimes (k+|j|)} 
\qquad 0 \le k+|j| \le m 
$$

We shall identify $b_0$ with an element in 
$A_{m/2} \subset A_m \subset L(X^{\otimes m})$.  
Apply Lemma \ref{lemma:T-epsilon} for  $m = n$ and 
$T = (b_0)^{1/2}$. Then  there exists
a positive element $a \in A$ such
$ \| \phi(\alpha^m(a))T\|^2 \ge \|T\|^2 -\varepsilon$
and $a \alpha^j(a) = 0$ for $j=1$, $\cdots$, $m$.
Define a positive operator $P = \alpha^m(a) \in A$.  Then 
$$
  \| P b_0 P\| = \| Pb_0^{1/2} \|^2 
                \ge \|b_0^{1/2}\|^2 -\varepsilon 
                = \| b_0 \| -\varepsilon
$$
For $j > 0$, we have   
\begin{align*}
 PS_xS_y^* P & =  \alpha^m(a)S_xS_y^* \alpha^m(a) 
             =  S_x \alpha^{m-(k+j)}(a) \alpha^{m-k}(a) S_y^* \\
             & =  S_x \alpha^{m-(k+j)}(a \alpha^j(a))S_y^* = 0 \\
\end{align*}
For $j<0$, we also have that $PS_xS_y^* P = 0$. 
Hence $Pb_jP = 0$ for $j \not= 0$.  

\end{proof}

\begin{thm}
Let $R$ be a rational function with $\deg R \geq 2$. 
Then the $C^*$-algebra ${\mathcal O}_R$ associated with
$R$ is simple and purely infinite.
  
\end{thm}

\begin{proof}
Let $w \in {\mathcal O}_R$ be any non-zero positive element.    
We shall show that there exist $z_1$, 
$z_2 \in {\mathcal O}_R $ such that $z_1^*w z_2 =I$.  
We may assume that $\|w\|=1$.
Let $E : {\mathcal O}_R \rightarrow  {\mathcal O}_R^{\gamma}$
be the canonical conditional expectation onto the fixed point 
algebra by the gauge action $\gamma$. 
Since $E$ is faithful, $E(w) \not= 0$.  
Choose  $\varepsilon$ such that 
\[
0 < \varepsilon < \frac{\|E(w)\|}{4} \  \text{ and } \ 
\varepsilon \|E(w) -3\varepsilon \|^{-1} \leq 1 .  
\]

There exists an element $c \in {\mathcal O}_X^{alg}$  
such that $ \|w - c^* c\| < \varepsilon$ and  
$\|c\| \le 1$.  Let $b = c^*c$. Then $b$ is decomposed  
as a finite sum $b = \sum_j b_j$ with 
$\gamma_t(b_j) =e^{ijt}b_j$.  
Since  $\|b\| \le 1$,  $\|b_0\| = \|E(b)\| \le 1$. 
By Lemma \ref{lemma:I-free}, there exists $P \in A$ with 
$0\le P \le I$ satisfying $Pb_jP = 0$ \qquad $(j\ne 0)$
and $\|Pb_0P\| \ge \|b_0\| -\varepsilon $. 
Then we have 

\begin{align*}
\| Pb_0 P \|  & \ge \|b_0 \| -\varepsilon 
              = \|E(b)\| -\varepsilon \\
              &  \ge \|E(w)\| -\|E(w) - E(b)\| -\varepsilon 
              \ge \|E(w)\| -2 \varepsilon 
\end{align*}
For $T := Pb_0P \in L(X^{\otimes m})$,  
there exists $f \in X^{\otimes m}$ with $\|f\|=1$ such that  
$$
 \|T^{1/2}f \|_2^2 = \|(f|Tf)_A\|  \ge \|T\| -\varepsilon 
$$
Hence we have 
$\|T^{1/2}f \|_2^2 \ge \|E(w)\| - 3 \varepsilon $.
Define $a = S_f^* T S_f = (f|Tf)_A \in A$.  
Then $\|a\| \ge \|E(w)\| -3 \varepsilon  > \varepsilon$. 
By Lemma \ref{lemma:u-epsilon},  there exists
$n \in \mathbb{N}$ and $u \in X^{\otimes n}$ 
sucn that  
$$
 \| u \| _2 \le (\|a\| - \varepsilon)^{-1/2} \qquad \text{and} 
\quad  S_u^* a S_u = I
$$
Then $\|u\| \le (\|E(w)\| -3 \varepsilon)^{-1/2}$. 
Moreover we have 

\begin{align*}
 \| S_f^* PwP S_f -a \| & = \|S_f^* PwPS_f -S_f^* T S_f\| \\
                        & = \|S_f^* PwPS_f -S_f^* Pb_0PS_f \| \\
                        & = \|S_f^* PwPS_f -S_f^* PbPS_f\| \\
                        & \le \|S_f\|^2 \|P\|^2 \|w -b\| <\varepsilon 
\end{align*}
Therefore   
\begin{align*}
   \|S_u^*S_f^* PwPS_fS_u -I\|  
  & = \|S_u^* S_f^* PwP S_fS_u - S_u^* a S_u \|  \\
  & \le \|u\|^2 \|S_f PwPS_f -a \| \\
  & < \|u\|^2 \varepsilon 
  \le \varepsilon \|E(w) -3\varepsilon \|^{-1} \le 1.
\end{align*}
Hence  
$S_u^*S_f^*PwPS_fS_u$ is invertible. Thus there  exists 
$v \in {\mathcal O}_X$ with $S_u^*S_f^*PwPS_fS_u v =I$. 
Put $z_1=S_u^*S_f^*P$ and $z_2=PS_fS_u v$. Then  
$z_1 w z_2 = I$.  

\end{proof}

\noindent
{\bf Remark}. J. Schweizer showed a nice criterion of the simplicity
of Cuntz-Pimsner algebras in \cite{S}: If a Hilbert bimodule $X$ is 
minimal and non-periodic, then ${\mathcal O}_X$ is simple.  Any 
$X$-invariant ideal $J$ of $A$  corresponds to a closed subset 
$K$ of $J_R$ with $R^{-1}(K) \subset K$. If $\deg R \geq 2$, then 
for  any $z \in J_R$ the backward orbit of $z$ is dense in $J_R$ 
by \cite{B} Theorem 4.2.7.  Therefore any $X$-invariant ideal $J$ of $A$ 
is $A$ or $0$, that is, $X$ is minimal.  Since $A$ is commutative 
and $L(X_A)$ is non-commutative, $X$ is non-periodic.  Thus 
Schweizer's theorem also implies that ${\mathcal O}_R$ is simple. 
Our theorem gives simplicity and pure infiniteness with a 
direct proof.

\begin{prop}
Let $R$ be a rational function with $\deg R \geq 2$. 
Then the $C^*$-algebra ${\mathcal O}_R$ associated with
$R$ is separable and nuclear, and satisfies the Universal  
Coefficient Theorem.  
\end{prop}

\begin{proof}Since ${\mathcal J}_X$ and  ${\mathcal T}_X$ are 
KK-equivalent to abelian $C^*$-algebras $I_X$ and $A$, 
the quotient ${\mathcal O}_X \cong {\mathcal T}_X/{\mathcal J}_X$ 
satisfies the UCT. ${\mathcal O}_X$ is also shown to be nuclear 
as in an argument of \cite{DS}.   
\end{proof}

\noindent
{\bf Remark}. If $\deg R \geq 2$, 
then the isomorphisms class of ${\mathcal O}_R$ is completely 
determined by the $K$-theory together with the class of the 
unit by the classification theorem by Kirchberg-Phillips \cite{Ki},
\cite{Ph}.

\section{examples}
If rational functions $R_1$ and $R_2$ are topologically conjugate, 
then their $C^*$-algebras  ${\mathcal O}_{R_1}$ and 
${\mathcal O}_{R_2}$ are isomorphic.  Therefore the K-groups 
$K_i({\mathcal O}_{R_1})$ and $K_i({\mathcal O}_{R_2})$ are
isomorphic.  Similarly 
$K_i({\mathcal O}_{R_1}(\hat{\mathbb C}))$ and \\
 $K_i({\mathcal O}_{R_2}(\hat{\mathbb C}))$ are
isomorphic. Moreover 
if rational functions $R_1$ and $R_2$ are topologically conjugate, 
then the gauge actions  are also conjugate.  
Therefore the K-groups of the fixed point algebras are also 
topologically conjugate invariant.
We investigate what kind of information of complex 
dynamical systems is described by the K-theory. 

We calculate the K-groups  by the following six-term exact 
sequence due to Pimsner \cite{Pi}.
$$
\begin{CD}
   K_0(I_X) @>{id - [X]}>> K_0(A) @>i_*>> K_0({\mathcal O}_R) \\
   @A{\delta _1}AA
    @.
     @VV{\delta _0}V \\
   K_1({\mathcal O}_R) @<<i_*< K_1(A) @<<{id - [X]}< K_1(I_X)
\end{CD}
$$

\noindent
{\bf Example 4.1.} Let $P(z) = z^n$ for $n \geq 2$.  
Then the Julia set $J_P$ is the unit circle $S^1$.  
The map $\phi : A \rightarrow L(X)$ can be identified 
with the $n$-times around embedding.
Hence we have 
$K_0({\mathcal O}_P) = {\mathbb Z} 
\oplus {\mathbb Z}/(n-1){\mathbb Z}$ 
and $K_1({\mathcal O}_P) = {\mathbb Z}$.
The fixed point algebra ${\mathcal O}_P^{\gamma}$ by the gauge 
action $\gamma$ is a Bunce-Deddence algebra of type $n^{\infty}$.

\medskip

\noindent
{\bf Example 4.2.} Let $P(z) = z^2 -2$. Then the Julia set $J_P$ is 
the interval $[-2,2]$ and  it contains a critical 
point 0. Since  $I_X = \{a \in C([-2,2]) ; a(0) = 0 \}$, 
$K_0(I_X) = K_1(I_X) = 0$.  Applying the above six term exact 
sequence with $K_0(A) = {\mathbb Z}$ and $K_1(A) = 0$,  we have 
$K_0({\mathcal O}_P) = {\mathbb Z}$ and 
$K_1({\mathcal O}_P) = 0$.  Since the identity $I$ of 
${\mathcal O}_P$ represents the generator of 
$K_0({\mathcal O}_P) = {\mathbb Z}$, 
the algebra ${\mathcal O}_P$ is isomorphic to
the Cuntz algebra ${\mathcal O}_{\infty}$.
$(J_P, P)$ is topologically conjugate to a tent map 
$([0,1], h)$ defined by 
$$
h(x) = \begin{cases}
 2x, & \qquad 0 \leq x \leq \frac{1}{2}, \\
-2x + 2, & \qquad \frac{1}{2} \leq x \leq 1.
\end{cases}
$$
Then $\frac{1}{2}$ is the only branched point of $h$. 
Deaconu and Muhly \cite{DM} associate a $C^*$-algebra 
$C^*([0,1],h)$ 
to the branched covering map $h: [0,1] \rightarrow [0,1]$. 
They exclude the branched point to construct their groupoid 
and the groupoid $C^*$-algebra $C^*([0,1],h)$ is not simple.
Moreover $K_0(C^*([0,1],h)) = {\mathbb Z}^2$ and 
$K_1(C^*([0,1],h)) = 0$.  Our $C^*$-algebra ${\mathcal O}_P$ 
is simple and purely infinite and is not isomorphic to their 
$C^*$-algebra $C^*([0,1],h)$.  The K-groups are different.

\medskip

\noindent
{\bf Example 4.3.} (quadratic polynomial) Let $P_c(z) = z^2 + c$. 
If $c$ is not in the Mandelbrot set ${\mathcal M} := 
\{c \in {\mathbb C} ; P_c^n(0) \text{ is bounded }\}$,
then ${\mathcal O}_{P_c}$ is isomorphic to the Cuntz algebra 
${\mathcal O}_2$.  In fact $(J_{P_c}, P_c)$ is topologically 
conjugate to the full two shift.  

If $c$ is in the interior of the main cardioid 
$C = \{ \frac{z}{2} - \frac{z^2}{4} \in {\mathbb C} ; |z| < 1 \}$, 
then the Julia set is homeomorphic to the unit circle $S^1$ 
(see, for example,  \cite{F}, page 211-212),  
and 
$P_c$ is topologically conjugate to $h(z) = z^2$ on $S^1$. Hence 
We have 
$K_0({\mathcal O}_P) = {\mathbb Z}$ and 
$K_1({\mathcal O}_P) = {\mathbb Z}$. 

\medskip

\noindent
{\bf Example 4.4.} Let $R$ be a rational function of degree $d \geq 2$.
Let $z_0$ be a (super)attracting fixed point of $R$. If all of 
the critical points of $R$ lie in the immediate attracting basin 
of $z_0$, then ${\mathcal O} _R$ is isomorphic to the Cuntz algebra 
${\mathcal O}_d$. In fact $(J_R, R)$ is topologically 
conjugate to the full d-shift by \cite{B} Theorem 9.8.1. See 
section 4 in \cite{PWY}.
For example, 
If $R(z) = \frac{2z^2 - 1}{z}$, then 
${\mathcal O}_R \cong {\mathcal O}_2$. 

\medskip

\noindent
{\bf Example 4.5.}  Tchebychev polynomials $T_n$ are defined by 
$\cos nz = T_n(\cos z)$. For example, $T_1(z) = z$, 
$T_2(z) = 2z^2 -1$.    Then the Julia set $J_{T_n}$ is 
the interval $[-1,1]$ and  $J_{T_n}$ contains  $n-1$ critical 
points for $n \geq 2$. Since $K_0(A) = {\mathbb Z}$, 
$K_1(A) = 0$, $K_0(I_X) = 0$ and $K_1(I_X) = {\mathbb Z}^{n-2}$,  
we have 
$K_0({\mathcal O}_{T_n}) = {\mathbb Z}^{n-1}$ and 
$K_1({\mathcal O}_{T_n}) = 0$.  Recall that,  
if the Julia set of 
a polynomial $P$ of degree $n \geq 2$ is the interval $[-1,1]$, 
then $P = T_n$ or $P = -T_n$ (\cite{B}, page 11). 

\medskip

\noindent
{\bf Example 4.6.} We consider a rational function $R$ by Lattes 
such that the Julia set is the entire Riemann sphere. Let 
$R(z) = \frac{(z^2 + 1)^2}{4z(z^2 - 1)}$.  Then  
$J_R = \hat{\mathbb C}$ contains six critical points 
and $K_0(I_X) = {\mathbb Z}$ and  $K_1(I_X) = {\mathbb Z}^5$. 
Therefore we have the exact sequence:
$$
\begin{CD}
    {\mathbb Z}@>{id - [X]}>> {\mathbb Z}^2 
@>i_*>> K_0({\mathcal O}_R) \\
   @A{\delta _1}AA
    @.
     @VV{\delta _0}V \\
   K_1({\mathcal O}_R) @<<i_*< 0 @<<{id - [X]}< {\mathbb Z}^5
\end{CD}
$$

Deaconu and Muhly \cite{DM} already obtained the similar diagram 
for their $C^*$-algebra $C^*(\hat{\mathbb C},R)$:

$$
\begin{CD}
    {\mathbb Z}@>{id - [X]}>> {\mathbb Z}^2 
@>i_*>> K_0(C^*(\hat{\mathbb C},R)) \\
   @A{\delta _1}AA
    @.
     @VV{\delta _0}V \\
   K_1(C^*(\hat{\mathbb C},R)) 
@<<i_*< 0 @<<{id - [X]}< {\mathbb Z}^8
\end{CD}
$$
Since they exclude branched points, the rank of  $K_0$-groups of them 
are different with ours and their $C^*$-algebra is also 
different with ours. 

\medskip

\noindent
{\bf Example 4.7.} Ushiki \cite{U} discovered a rational 
function whose Julia set is homeomorphic to the  
Sierpinski gasket.  See also \cite{Kam}.  For example,   
let $R(z) = \frac{z^3-\frac{16}{27}}{z}$. Then $J_R$ is 
homeomorphic to the  Sierpinski gasket K and $J_R$ contains 
three critical points.  Recall that the usual Sierpinski 
gasket $K$ is constructed by three contractions 
$\gamma _1, \gamma _2, \gamma _3$ on the regular triangle $T$
in ${\mathbb R}^2$ with three vertices $P = (1/2,\sqrt{3}/2)$, 
$Q = (0,0)$ and $R = (1,0)$ such that 
$\gamma _1(x,y) = (\frac{x}{2} + \frac{1}{4}, \frac{y}{2} + \frac{\sqrt{3}}{4})$, $\gamma _2(x,y) = (\frac{x}{2}, \frac{y}{2})$, 
$\gamma _3(x,y) = (\frac{x}{2} + \frac{1}{2}, \frac{y}{2})$. 
Then a self-similar set $K \subset T$ satisfying 
$K = \gamma _1(K) \cup \gamma _2(K) \cup \gamma _3(K)$ is called 
a Sierpinski gasket.  But these three contractions are not 
inverse branches of a map, because 
$\gamma _1(Q) = \gamma _2(P)$. 
Therefore we need to modify the 
construction of contractions.  Put  
$\tilde{\gamma _1} = \gamma _1$, 
$\tilde{\gamma _2} = \alpha _{-\frac{2\pi}{3}} \circ \gamma _2$,  
and $\tilde{\gamma _3} = \alpha _{\frac{2\pi}{3}} \circ \gamma _3$, 
where $\alpha _{\theta}$ is a rotation by the angle $\theta$. 
Then $\tilde{\gamma _1}, \tilde{\gamma _2},  \tilde{\gamma _3}$ 
are inverse branches of a map $h: K \rightarrow K$, which is 
conjugate to $R : J_R \rightarrow J_R$.  Two correspondences 
$C = \cup _i \{(\gamma _i(z),z) ; z \in K \}$ and 
$\tilde{C} = \cup _i \{(\tilde{\gamma _i}(z),z) ; z \in K \}$
generate Hilbert bimodule $Z$ and $\tilde{Z}$ over $C(K)$. 
Then $C^*$-algebra ${\mathcal O}_R \cong {\mathcal O}_{\tilde{Z}}$
and $K_0({\mathcal O}_R)$ contains a torsion free element.  But 
$C^*$-algebra ${\mathcal O}_Z$ is isomorphic to the Cuntz algebra 
${\mathcal O}_3$.  Therefore  $C^*$-algebra ${\mathcal O}_Z$ and 
${\mathcal O}_R$ are not isomorphic.  See \cite{KW}.

\section{Fatou set and the corresponding ideal}

Recall that the Fatou set $F_R$ of a rational functon $R$ 
is the maximal open subset of the Riemann sphere 
$\hat{\mathbb C}$ on which $(R^n)_n$ is equicontinuous, 
and the Julia set $J_R$ of $R$ is the 
complement of the Fatou set in $\hat{\mathbb C}$. We consider 
the corresponding decomposition for the $C^*$-algebra 
${\mathcal O}_R(\hat{\mathbb C})$, which was first pointed out 
by Deaconu and Muhly \cite{DM} in the case of their construction.
Let $B = C(\hat{\mathbb C})$ and consider the ideal 
$I = \{b \in B ; b|J_R = 0\}$ of $B$, so  $I \cong C_0(F_R)$ 
and $B/I \cong C(J_R)$.  We consider a submodule and quotient module 
of a Hilbert bimodule $Y= C(\graph R)$ over $B$.  
The right Hilbert $I$-module 
$YI:= \{ fb  \in Y ; f \in Y, \ b \in I \}$ is also described as 
$Y_I : = \{f \in Y ; (f|g)_B \in I \text{ for all } g \in Y \}$. Since 
$(f|f)_B \in I$ means that 
$$
(f|f)_B(y) = \sum _{x \in R^{-1}(y)} e(x) |f(x,y)|^2 = 0  
$$
for all $y \in J_R$, we have 
$$
YI = \{f \in Y ; f(x,y) = 0  \text{ for all } (x,y) \in \graph R|_{J_R} \}. 
$$
because $J_R$ is complete invariant.  

Invariant ideals for bimodules are introduced in \cite{KPW1} by Pinzari and 
ours and 
developed by Fowler, Muhly and Raeburn \cite{FMR} in general case.
In our situation   
$I$ is a $Y$-invariant ideal of $B$, i.e., 
$\phi (I)Y \subset YI$.  In fact, the condition is equivalent to 
that $(f|\phi (a)g)_B \subset I$ for any $a \in I$ and $f,g\in Y$, 
and it is easily checked as 
$$
(f|\phi (a)g)_B(y) 
= \sum _{x \in R^{-1}(y)} e(x) \overline{f(x,y)}a(x)g(x,y) = 0
$$
for $y \in J_R$, because $x \in R^{-1}(J_R) = J_R$ and so 
$a(x) = 0$ for $a \in I$. 

Therefore $Y/YI$ is naturally a Hilbert bimodule over 
$B/I \cong C(J_R)$.  We can identify  $Y/YI$ with a 
bimodule $X = C(\graph R|_{J_R})$ over $A = C(J_R)$. 
  
\begin{thm}
Let $R$ be a rational function,
$B = C(\hat{\mathbb C})$, $A = C(J_R)$ 
and  $I = \{b \in B ; b|J_R = 0\}$.
Then the ideal ${\mathcal I} (I)$ generated by $I$ 
in ${\mathcal O}_R(\hat{\mathbb C})$ is Morita equivalent to 
${\mathcal O}_R(F_R)$ and the quotient algebra 
${\mathcal O}_R(\hat{\mathbb C})/{\mathcal I} (I)$ is 
canonically isomorphic to ${\mathcal O}_R = {\mathcal O}_R(J_R)$.
\end{thm}

\begin{proof} We apply a result by Fowler, Muhly and Raeburn 
\cite{FMR} Corollary 3.3. But we have to be careful, because 
$\phi (B)$ is {\it not} included in $K(Y)$. The only thing 
we know is that   
the quotient algebra 
${\mathcal O}_R(\hat{\mathbb C})/{\mathcal I} (I)$ is canonically 
isomorphic to the relative Cuntz-Pimsner algebra 
${\mathcal O}(q^I(I_Y), Y/YI)$.  Here $q^I : B \rightarrow B/I \cong A$
is the canonical quotient map. By Proposition \ref{prop:critical},  
$I_Y = \{b \in B ; b \ \text{vanishes on } C \}$.
and  $I_X = \{a \in A ; a \ \text{vanishes on } C \cap J_R \}$. 
Hence we can identify $q^I(I_Y)$ with $I_X$.  Since we also 
identify a bimodule $Y/YI$ over $B/I$ with a 
bimodule $X = C(\graph R|_{J_R})$ over $A = C(J_R)$, 
we have that ${\mathcal O}(q^I(I_Y), Y/YI)$ is isomorphic to 
${\mathcal O}_R(J_R)$.  
\end{proof}

\section{Lyubich measure and KMS state}

Lyubich \cite{L} constructed an invariant measure $\mu$ for 
a rational function $R$, called a Lyubich measure, 
whose support is the Julia set.  It derives a lower bound for 
the Hausdorff dimension of the Julia set. Lyubich showed 
that $\mu$ is a unique measure of maximal entropy and 
the entropy value $h_{\mu}(R)$ is equal to the topological entropy
$h(R)$ of $R$, which is $\log (\deg R)$. 
We show that the Lyubich measure gives a unique KMS state on 
the $C^*$-algebra ${\mathcal O}_R$ for the gauge action at 
inverse temperature $\log (\deg R)$ if the Julia set contains 
no critical points. 

Let $d = \deg R$ and $\delta _x$ be the Dirac measure for $x$ 
on the Riemann sphere.  For any $y \in \hat{\mathbb C}$ 
and each $n \in {\mathbb N}$, we define 
a probability measure $\mu _n^y$ by 
$$
\mu _n^y = \sum _{\{ x \in \hat{\mathbb C} ; R^n(x) = y \}} 
           d^{-n}e(x)  \delta _x \ . 
$$ 
Lyubich showed that the sequence $(\mu _n^y)_n$ converges 
weakly to a measure $\mu$. The Lyubich measure $\mu$ is independent 
of the choice of $y$, the support of $\mu$ is the Julia set
$J_R$ and $\mu (E) = \mu (R^{-1}(E))$ for any Borel set $E$. 
Hence $\int a(R(x)) d\mu (x) = \int a(x) d\mu (x)$. 

\begin{thm}
Let $R$ be a rational function with $\deg R \geq 2$. 
Consider the $C^*$-algebra ${\mathcal O}_R$ associated with
$R$ and the gauge action 
$\gamma :  {\mathbb R} \rightarrow Aut \ {\mathcal O}_R$
such that $\gamma_t(S_f) = e^{it}S_f$ for $f \in X$. Suppose 
that the Julia set contains no critical points, (in particular, 
suppose that $R$ is hyperbolic). 
Then $({\mathcal O}_R, \gamma)$ has a KMS state at the inverse 
temperature $\beta$ if and only if $\beta = \log (\deg R)$ and 
the corresponding KMS state $\varphi$ is also unique.  Moreover 
the restriction of the state $\varphi$ to $A = C(J_R)$ is given 
by the Lyubich measure $\mu$ such that $\varphi (a) = \int a d\mu$
for $a \in A$.
\end{thm}

\begin{proof}
Since the Julia set contains no critical points, $A$-module 
$X$ has a finite right basis $\{u_1, \dots , u_n\}$.  Let 
$\varphi$ be a state on ${\mathcal O}_R$ and $\varphi _0$ the 
restriction of $\varphi$ to $A$.  As in \cite{PWY} section 3, 
$\varphi$ is a KMS state at the inverse temperature $\beta$ 
if and only if $\varphi _0$ satisfies that 
$$
(*) \ \ \ \ \  \varphi _0(\sum _i (u_i|au_i)_A) = 
  e^{\beta}\varphi _0 (a)    \ \ \ \text{for all } a \in A.
$$
Moreover any state $\varphi _0$ on $A$ satisfying (*) extends 
uniquely to a KMS state $\varphi$ at the inverse temperature 
$\beta$. In fact, define a state $\varphi _n$ on the $C^*$-algebra
$B_n $ generated by 
\[
\{S_{x_1} \dots S_{x_n}S_{y_n}^* \dots S_{y_1}^* ;  
x_1, \dots x_n,  y_1, \dots y_n \in X \},
\] 
by the induction: 
$$
\varphi _{n+1}(T) =  e^{-\beta} \varphi _n (\sum _i S_{u_i}^*TS_{u_i})
$$
for $T \in B_n$.  The sequence $(\varphi _n)_n$ gives a state 
$\varphi _{\infty}$ on the fixed point algebra 
${\mathcal O}_R^{\gamma}$. Let 
$E : {\mathcal O}_R \rightarrow  {\mathcal O}_R^{\gamma}$
be the canonical conditional expectation.  Then the KMS state 
$\varphi$ is defined by $\varphi = \varphi _{\infty} \circ E.$ 

Define a complete positive map $h : A \rightarrow A$ by 
$h(a) = \sum _i (u_i|au_i)_A$ for $a \in A$.  Then we have 
$$
h(a)(y) = \sum _{x \in R^{-1}(y)} a(x)
$$
for $y \in J_R$.  In fact, since $\{u_1, \dots , u_n\}$ is a 
finite basis for $X$, $\sum _i u_i(u_i|f)_A = f$ for any 
$f \in X$.  Therefore for $y = R(x)$, 
$$
\sum _i u_i(x,y) \sum _{x' \in R^{-1}(y)} 
\overline{u_i(x',y)}f(x',y) = f(x,y).
$$
This implies that $\sum _i u_i(x,y) \overline{u_i(x',y)} = 
\delta _{x,x'}$. Hence we have 
$$
\sum _i u_i(x,y) \overline{u_i(x,y)} = 1 .
$$
Therefore 
\begin{align*}
h(a)(y) & = \sum _i (u_i|au_i)_A(y) \\ 
& = \sum _i \sum _{x \in R^{-1}(y)} 
\overline{u_i(x,y)} a(x)u_i(x,y) = \sum _{x \in R^{-1}(y)}a(x).
\end{align*}
Hence the condition (*) is written as 
$$
\varphi _0(h(a)) = e^{\beta} \varphi _0(a) 
\ \ \ \text{for all } a \in A.
$$
Putting $a = 1$, we have that  $\deg R = e^{\beta}$, that is, 
$\beta$ need to be $\log (\deg R)$.  For any $a \in A$, 
$(e^{-\beta}h)^n(a)$ converges uniformly to a constant 
$\int a d\mu$ by \cite{L}.  Thus 
$\int e^{-\beta}h(a) d\mu = \int a d\mu$.  
Hence Lyubich measure $\mu$ gives a KMS state  
at the inverse temperature $\log \deg R$. Take 
another KMS state $\tau$.  Then  
$\tau ((e^{-\beta}h)^n(a) = \tau (a)$ converges to $\int a d\mu$.
Thus $\tau(a) = \int a d\mu$. It shows that KMS state is unique. 

\end{proof}


%

\end{document}